\documentclass{amsart}

\usepackage{amsmath,amsthm,amscd,amssymb,latexsym,verbatim}

\newcommand{\R}{\ensuremath{\mathbb{R}}}   
\newcommand{\C}{\ensuremath{\mathbb{C}}}   
\newcommand{\T}{\ensuremath{\mathbb{T}}}   
\newcommand{\Q}{\ensuremath{\mathbb{Q}}}   

\renewcommand{\S}{\ensuremath{\mathbb{S}}}

\theoremstyle{plain}
\newtheorem{theorem}{Theorem}[section]

\theoremstyle{definition}

\theoremstyle{remark}
\newtheorem{remark}[theorem]{Remark}

\numberwithin{equation}{section} \allowdisplaybreaks

\def \d {\delta}
\def \g {\gamma}
\def \e {\varepsilon}
\def \f {\varphi}
\def \F {{\bf \Phi}}

\def \l {\lambda}

\def \n {\nabla}
\def \s {\sigma}

\def \O {\Omega}

\def \lan {\langle}
\def \ran {\rangle}
\def \p {\partial}
\def \ra {\rightarrow}
\def \ss {\subset}

\newcommand{\der}[2]{(#1 \cdot \nabla) #2}

\DeclareMathOperator{\re}{Re} %
\DeclareMathOperator{\rg}{Rg} %
\DeclareMathOperator{\Ker}{Ker} %

\def \sol { L^2_{\mathrm{div}} }

\def \ess {\mathrm{ess}}

\begin{document}

\title{Continuous spectrum of the 3D Euler equation is a solid annulus}

 \author{R. Shvydkoy}
\address{Department of Mathematics, Statistics, and Computer Science, \\
851 S Morgan St., M/C 249,\\
University of Illinois at Chicago, Chicago, IL 60607}
\email{shvydkoy@math.uic.edu}
\thanks{The work is partially supported by the US National Science
Foundation Grant DMS-0907812.}

\begin{abstract}

In this note we give a description of the continuous spectrum of the linearized Euler equations in three dimensions. Namely, for all but countably many times $t\in \R$, the continuous spectrum of the evolution operator $G_t$ is given by a solid annulus with radii $e^{t\mu}$ and $e^{t M}$, where $\mu$ and $M$ are the smallest and largest, respectively, Lyapunov exponents of the corresponding bicharacteristic-amplitude system of ODEs.

{\bf Le spectre continu de l'equation d'Euler 3D est un anneau}.
Nous donnons dans cette article une description du spectre continu
de l'equation d'Euler linearisee en dimension 3.
Precisement, pour presque tout $t\in \R$, le spectre continu de l'operateur
d'evolution  $G_t$ est constitue d'un anneau de rayons $e^{t\mu}$ et $e^{t M}$, ou $\mu$ et $M$ sont respectivement le plus petit et le plus grand
exposant de Lyapunov du systeme de EDO bicharacteristique-amplitude associe.

\end{abstract}

\keywords{Euler equation, Fredholm spectrum, essential spectrum, dynamical spectrum, cocycle}
\subjclass[2000]{47D03; 76E09 }

\maketitle

\section{Background}

This note addresses the geometric structure of the continuous spectrum (in Browder or Fredholm sense) of the linearized three-dimensional Euler equation. The question has its origin in the so-called elliptic instabilities of ideal fluids found implicated in transition to turbulence of some lamina flows in a subcritical range of Reynolds numbers (see Bayly \cite{Bayly}, Patera and Orszag \cite{PO}, Pierrehumbert \cite{P}). The works of Lipschitz and Hameri \cite{LH}, and  Friedlander and Vishik \cite{FV} laid the basis for a general approach to shortwave instabilities which incorporated elliptic and a broad range of other equilibria, including those possessing exponential stretching of trajectories. In \cite{V} Vishik draws a connection between shortwave instabilities and the essential spectrum of the linearized Euler equation, as opposed to the point spectrum, and describes the essential spectral radius in terms of the maximal Lyapunov-Oseledets exponent of a bicharacteristic-amplitude system of ODEs, BAS for short, (see below). Further investigation on the structure of the spectrum continued in a series of papers \cite{LV,SL,SL2,S}, where the Sacker-Sell theory of linear skew-product flows was introduced into the subject. It was shown that in general the shortwave instabilities and the corresponding points of the essential spectrum of the Euler equations are linked to the dynamical spectrum of the BAS (in the Sacker-Sell sense). At the same time the spectrum in two dimensions was completely described in \cite{SL} by a direct construction of approximate eigenfunctions. Thus, the spectrum of the Euler semigroup is a solid annulus centered at the origin, while the essential spectrum of the generator is the corresponding vertical band, so that the two are related by the Spectral Mapping Theorem.  Over the past several years attempts were made to find a similar description in three dimensions. In this note we offer a complete answer to this question on the level of the semigroup proving that, just as in two dimensions, its essential spectrum is given by a solid annulus. A partial description of the spectrum for the generator was also provided in \cite{S}. However, the spectral mapping theorem in this case remains open, leaving a possibility for further investigation.

\section{Technical description of the result}
Let $u_0 \in (C^2(\T^3))^3$ be a solution to the stationary Euler equations
\begin{equation}\label{e:ss}
\der{u_0}{u_0} + \n p = 0,\quad
\n \cdot u_0 = 0,
\end{equation}
on the three dimensional torus $\T^3$. Evolution of small perturbations to $u_0$ are described up to the linear approximation by the system
\begin{gather}\label{lee}
v_t = - \der{u_0}{v} - \der{v}{u_0} - \n q, \quad
\n \cdot v = 0 .
\end{gather}
The solution map for \eqref{lee} is given by a $C_0$-group $\{G_t\}_{t\in \R}$ of operators acting on the energy space $\sol = (\sol(\T^3))^3$ of divergence-free vector fields. The method of geometric optics projects evolution of  oscillating localized perturbations of the form
 $v(x,t) = b(x,t) e^{i S(x,t) /\d} + O(\d)$,  $ 0<\d \ll 1$,
where the initial amplitude $b(x,0)=b_0(x)$ and phase $S(x,0) = \xi_0 \cdot x$ are given, and $\xi_0 \cdot b_0(x) = 0$ for all $x\in\T^3$. If written in the Lagrangian coordinates
of the flow $u_0$, $x_0 \ra \f_t(x_0)$, the evolution of the vectors $b(t) = b(\f_t(x_0),t)$ and $\xi(t) = \n S(\f_t(x_0),t)$ is governed by the bicharacteristic-amplitude system (BAS) of ODEs:
\begin{subequations}\label{BAS}
\begin{align}
    x_t &= u_0(x), & \quad x(0) & = x_0, \label{BASx} \\
    \xi_t &= -\partial u_0(x)^\top \xi, & \quad \xi(0) & = \xi_0 \perp b_0\label{BASxi}\\
    b_t &=-\partial u_0(x)b + 2(\partial u_0(x)b, \xi)\xi |\xi|^{-2}, &\quad b(0) & = b_0. \label{BASb}
\end{align}
\end{subequations}
One can see that $b \cdot \xi = 0$ is a conservation law of \eqref{BAS}, which is consistent with the incompressibility of the flow. We view equations \eqref{BASx} - \eqref{BASxi} as a Hamiltonian system over the cotangent bundle $T^*\T^3 = \T^3 \times \R^3$, with the Hamiltonian $H(x,\xi) = u_0(x) \cdot \xi$. The corresponding flow projected onto the compact space $\O = \T^n \times \S^{n-1}$ is given by
\begin{equation}\label{chi-flow}
\chi_t(x_0,\xi_0) = \left( \f_t(x_0) , \frac{\p \f_t^{-\top}(x_0) \xi_0}{|\p \f_t^{-\top}(x_0) \xi_0|} \right).
\end{equation}
As the amplitude equation \eqref{BASb} is homogeneous in $\xi$ it can be written in the form $b_t = a_0(\chi_t(x_0,\xi_0)) b$. Therefore, the fundamental matrix solution of \eqref{BASb}, $B_t(x,\xi)$ defines a cocycle over the flow $\{\chi_t\}$ on the fibre bundle $\mathcal{F}$ with base $\O$ and fibers given by $F(x,\xi) = F(\xi) = \{ b\in \C^3: b \cdot \xi = 0\}$. We call it the $b$-cocycle. Let $\Sigma_B$ denote the dynamical spectrum of the $b$-cocycle,
that is, the set of $\l \in\R$ such that  the rescaled cocycle $B_t^\l = e^{-\l t}B_t$ {\em does not} have exponential dichotomy (uniform exponential hyperbolicity) on the bundle, see \cite{SS}.
In the resent paper \cite{SL2}, we showed a key connection between the Fredholm spectrum of $G_t$ and the dynamical spectrum $\Sigma_B$, expressed by the following two identities:
\begin{equation}\label{almostal}
\s_{\F}(G_t; \sol) = \{ z \in \C: |z| \in e^{t \Sigma_B} \},
\end{equation}
for all but countably many $t\in \R$, and
\begin{equation}\label{al}
| \s_{\F}(G_t; \sol) | = e^{t \Sigma_B},
\end{equation}
for all $t$ without exceptions. Here $\s_{\F}$ stands for the Fredholm spectrum over $\sol$.
The classical results of Sacker and Sell on the structure of a dynamical spectrum, see \cite{SS}, in our case reveal that  $\Sigma_B$ is either a single interval $[\mu, M]$ or the union of at most two disjoint intervals, as limited by the dimension of the fibers $F(\xi)$. The main contribution of this note is to add one last technical piece to \eqref{almostal} - \eqref{al}, which is to show that $\Sigma_B$ is in fact a single interval.  We thus obtain the following theorem.
\begin{theorem}\label{t:main} Let $\mu$ and $M$ be the minimal and maximal Lyapunov-Oseledets exponents of the $b$-cocycle, respectively. Then the following identity
\begin{equation}\label{e:annulus}
 \s_\ess(G_t;\sol) = \s_{\F}(G_t;\sol) = \{ z \in \C: e^{t\mu} \leq |z| \leq e^{t M} \}
\end{equation}
holds for all but countably many $t\in \R$, while the identity
\begin{equation}\label{e:hull}
    | \s_\ess(G_t;\sol) | = | \s_{\F}(G_t;\sol) | = [e^{t\mu}, e^{t M} ]
\end{equation}
holds for all $t\in \R$ without exception.
\end{theorem}
Here $\s_{\ess}$ stands for the essential spectrum in the Browder sense \cite{Browder}. We added it to the statement as this is the type of spectrum that was examined in all the previous works. As we see from \eqref{e:annulus} it is in fact the same as the Fredholm spectrum for almost all $t$. The identities \eqref{e:annulus} - \eqref{e:hull} for $\s_{\ess}$ follow immediately from the identities for $\s_{\F}$ due to the general inclusion $\s_{\F} \ss \s_{\ess}$ and the Nussbaum Theorem \cite{Nuss} on the equality between the spectral radii (applied also to the inverse of $G_t$).

\section{Connectedness of the dynamics spectrum}
Let us suppose, on the contrary, that $\Sigma = [\mu,\g_1] \cup [\g_2,M]$, where $\g_1<\g_2$.  Then there exists a continuous projection-valued function
\begin{equation}\label{e:proj}
P(x,\xi) :  F(\xi) \ra F(\xi)
\end{equation}
corresponding to the exponential splitting of the cocycle $B_t$. So, for a small $\e >0$ and $\l \in (\g_1,\g_2)$ there are constants $c,C>0$ such that
\begin{align}
\| B_t^\l(x,\xi) b \| & \leq C e^{-\e t}\|b\|, \, t >0,  \text{ for all } b \in \rg P(x,\xi); \label{exp1}\\
\| B_t^\l(x,\xi) b \| & \geq c e^{\e t} \|b\|, \, t>0,  \text{ for all } b \in \Ker P(x,\xi).\label{exp2}
\end{align}
Moreover, relations \eqref{exp1} and \eqref{exp2} characterize the range and kernel of $P(x,\xi)$, respectively.

Let us fix an arbitrary $x_0 \in \T^3$, and consider the splitting
\begin{equation}\label{e:split}
F(\xi) = \rg P(x_0,\xi) \oplus \Ker P(x_0,\xi).
\end{equation}
Since the dimension of $F(\xi)$ over $\C$ is two, and \eqref{e:split} is non-trivial, the dimension  of $\rg P(x_0,\xi)$ is one. Let us now fix a finite local coordinate chart $\{U_j\}_{j\in J}$ of $\S^2$. For each $j\in J$ we can find a continuous mapping
$$
\xi \ra b_\C^j(\xi) \in F(\xi), \quad \xi \in U_j,
$$
with $\| b_\C^j(\xi) \| = 1$,  such that
\begin{equation}\label{e:spanC}
\rg P(x_0,\xi) = [ b_\C^j(\xi) ]_\C,
\end{equation}
where $[\cdot]_\C$ denotes the linear span over $\C$.

We will now de-complexify the spaces given by \eqref{e:spanC}.  To this end, let us write
$$
b_\C^j (\xi)= b_{re}^j(\xi) + i b_{im}^j(\xi),
$$
where $b_{re}^j(\xi)$ and $b_{im}^j(\xi)$ are real vectors. Since $b_\C^j(\xi) \cdot \xi = 0$, and $\xi$ is real, we also have
$$
b_{re}^j(\xi), \, b_{im}^j(\xi) \in F(\xi).
$$
Furthermore, the cocycle $B_t$, being the fundamental matrix of solutions of a system with real coefficients, maps real vectors to real vectors. Thus, we obtain
$$
\|B_t^\l (x_0,\xi) b^j_{re}(\xi) \| + \|B_t^\l (x_0,\xi) b^j_{im}(\xi) \| \sim \|B_t^\l (x_0,\xi) b^j_\C(\xi) \|  \leq  C e^{-\e t},
$$
for all $t >0$.  This implies that both vectors $b_{re}^j(\xi)$ and $b_{im}^j(\xi)$ belong to $\rg P(x_0,\xi)$, span it over $\C$, and are linearly dependent over $\R$. Let us define the map
\begin{equation}\label{Lj}
L_j(\xi)  = [ b_{re}^j(\xi), b_{im}^j(\xi) ]_{\R},  \quad \xi \in U_j.
\end{equation}
We see that each $L_j(\xi)$ is a one-dimensional subspace of the real tangent plane to the surface of the sphere $\S^2$ at $\xi$. It is not hard to see that $L_j$'s agree on intersections of the charts.  Indeed, let $\xi \in U_j \cap U_k$, for some $j \neq k$.  According to \eqref{e:spanC}, vectors $b_\C^j(\xi)$ and $b_\C^k(\xi)$ are linearly dependent over $\C$. Thus, for some $z = x+iy$,  $b_\C^j (\xi) =z b_\C^k(\xi)$, which implies linear dependencies
\begin{align}
b_{re}^j(\xi) &= x b_{re}^k(\xi) -  y b_{im}^k(\xi), \\
b_{im}^j(\xi) & = y b_{re}^k(\xi)+xb_{im}^k(\xi).
\end{align}
It is then clear that $L_j(\xi) = L_k(\xi)$.  As a result, \eqref{Lj} defines a continuous one-dimensional subbundle of the tangent bundle of $\S^2$. According to \cite[Theorem 27.16]{Steenrod}, one can then select a continuous non-vanishing tangent field on the sphere, which contradicts the classical Hairy Ball Theorem.

\begin{remark}
A simple example disclosed in \cite{SL2} shows that the "all but countable many" condition in Theorem~\ref{t:main} is essential. Consider a two dimensional parallel shear flow with constant profile, $u_0 = \lan U,0 \ran$. The linearized equation \eqref{lee} takes the form $v_t = - U \p_{x_1} v$. We then obtain $G_t v_0(x_1,x_2) = v_0(x_1 - U t \mod 2\pi, x_2)$. So, the spectrum of $G_t$ is the unit circle for all $t \not \in \pi U^{-1} \Q$, and a finite set of the circle otherwise. 
\end{remark}

\begin{remark}It is desired to obtain a similar description of the essential spectrum for the generator $L$, which would exponentially correspond to the solid annulus. So far partial results has been given in \cite{S}. Let us consider the closed subset $\O^\perp = \{(x,\xi) \in \O: u_0(x) \cdot \xi = 0\}$. The flow $\chi_t$ leaves $\O^\perp$ invariant. We can therefore consider the restrictions of the flow, $\chi_t^\perp$, and the cocycle $B_t^\perp$ on the fiber bundle $\mathcal{F}^\perp$ on $\O^\perp$ with the same fibers. Let $\Sigma_{B^\perp}$ be the corresponding dynamical spectrum of the cocycle $\{B_t^\perp\}$. It was shown in \cite{S} that $\Sigma_{B^\perp} \ss \re( \s_{\ess}(L))$.  We can see now that if there is a stagnation point $x_0\in\T^3$ of the flow, i.e. $u_0(x_0)=0$, then the corresponding $x_0$-slice of $\O^\perp$ is (again) the sphere $\S^2$. So, the argument above applies with this specific choice of $x_0$, which implies that $\Sigma_{B^\perp}= [\mu^\perp,M^\perp]$, for some $\mu^\perp \geq \mu$ and $M^\perp \leq M$. It remain to be seen, however, whether $\Sigma_{B^\perp}$ is in fact the same as $\Sigma_{B}$. So far, it has been established only in the two dimensional case, see \cite{LV,SL}.
\end{remark}


\end{document}